\documentclass[12pt]{amsart}

\usepackage{amsthm}
\usepackage{comment}

\usepackage{geometry}
\usepackage{amsmath,amssymb,mathrsfs}
\usepackage{graphicx}
\usepackage{enumerate}
\usepackage{enumitem}
\usepackage{mathtools}
\usepackage{xcolor}

\usepackage{hyperref}

\numberwithin{equation}{section}

\newtheorem{thm}{Theorem}[section]
\newtheorem{lem}[thm]{Lemma}
\newtheorem{prop}[thm]{Proposition}

\newtheorem{definition}[thm]{Definition}
\newtheorem{rmk}[thm]{Remark}
\newtheorem{rmks}[thm]{Remarks}

\newcommand\Var{\mathrm{Var}}

\newcommand{\E}{{\mathbb E}}
\newcommand{\T}{{\mathbb T}}
\newcommand{\N}{{\mathbb N}}
\newcommand{\R}{{\mathbb R}}

\author{Jacques Benatar}
\author{Alon Nishry}
\address{ Department of Mathematics and Statistics, University of Helsinki,
PO Box 68, FI-00014 Helsinki, Finland }
\email{jacques.benatar@helsinki.fi}
\address{School of Mathematical Sciences,
Tel Aviv University,
Tel Aviv 69978,
Israel}
\email{alonish@tauex.tau.ac.il}

\begin{document}

\begin{abstract}
For $X(n)$ a Steinhaus random multiplicative function, we study the maximal size of the random Dirichlet polynomial
\[
D_N(t) = \frac1{\sqrt{N}} \sum_{n \leq N} X(n) n^{it},
\]
with $t$ in various ranges. In particular, for fixed $C>0$ and any small $\varepsilon>0$ we show that, with high probability, 
\[
 \exp( (\log N)^{1/2-\varepsilon} ) \ll \sup_{|t| \leq N^C} |D_N(t)| \ll \exp(  (\log N)^{1/2+\varepsilon}).
\]
\end{abstract}

\title[Dirichlet polynomials with random multiplicative coefficients]{Extremal bounds for Dirichlet polynomials with random multiplicative coefficients}

\maketitle

\section{Introduction}

\subsection{Set-up and the main result}

Our central object of study is the normalised random Dirichlet polynomial 
\begin{equation}\label{Ddef}
D_N(t) = \frac1{\sqrt{N}} \sum_{n \leq N} X(n) n^{it},
\end{equation}
generated by coefficients $(X(n))_{n \in \N}$ which form a Steinhaus \textsl{random multiplicative function}, or RMF for short. We recall their construction: letting $(X(p))_{p}$ be a sequence  of i.i.d. random variables, indexed over the primes and uniformly distributed on the unit circle, we set
\begin{equation}\label{Xndef}
X(n)=\prod_{p^e || n} X(p)^e
\end{equation}
for each natural number $n \geq 1$. Here $p^e$ is the largest power of $p$ dividing $n$. With this definition, $X(n)$ forms a completely multiplicative sequence of \textsl{dependent} variables.

As a complement to the work of Rodgers and the authors \cite{JAB}, in which the distribution of the trigonometric polynomial with coefficients $X(n)$ was investigated, the purpose of this note is to study the large values of $|D_N(t)|$, with $t$ in various ranges.\\ 
Before stating our main result, we recall that a sequence of events $E_n$ is said to occur asymptotically almost surely if $\mathbb{P} (E_n) = 1-o(1)$ as $n \rightarrow \infty$. It will also be convenient to write $\log_k $ for the $k$-fold iterated logarithm.
\begin{thm}\label{mainprop}
Suppose that $C: (3,\infty) \rightarrow (0,\infty)$ satisfies the growth conditions
\begin{equation}\label{Cgrowth}
 \frac{(\log_2 x)^{9}}{\log x}  \leq C(x) \leq (\log x)^{\gamma}
\end{equation}
for some fixed exponent $0< \gamma < 1$ and all sufficiently large $x$. Consider the supremum
$$\mathcal{S}(N,C)=\sup_{|t| \leq N^{C(N)} } |D_N(t)|,$$
where $D_N$ is the random Steinhaus Dirichlet polynomial defined in \eqref{Ddef}. Then for any fixed $\varepsilon >0$, the estimates 
\begin{equation}\label{mainestimate}
 \exp \left( \frac{B \sqrt{C(N)  \log N} }{(\log_2 N)^{2} } \right) \ll_{ \gamma} \mathcal{S}(N,C) \ll \exp \left( (\tfrac{3}{2} +\varepsilon ) \sqrt{C(N)  \log N \log_2 N} \right)
\end{equation}
hold asymptotically almost surely as $N \rightarrow \infty$, with $B > 0$ an absolute constant. The upper bound holds uniformly over all $C>1$, in the sense that there is no restriction on the size of $\gamma \in (0,\infty)$.
\end{thm}

\begin{rmks} $\,$

\begin{enumerate}[label=(\alph*)]
\item We will see in section 3 below that the upper bound in Theorem~\ref{mainprop} follows from a basic moment estimate for $D_N(0)$. As such, it is certainly not a new result (e.g. Granville and Soundararajan \cite[Theorem 4.1]{GS}). A far more elaborate and delicate treatment of the moments of $D_N(0)$ can be found in 
Harper \cite{HarII}. Employing the moment bound from \cite{HarII} would likely yield a small improvement of the constant $3/2$ in Theorem ~\ref{mainprop}, but we do not pursue the matter here. 
\item The bounds \eqref{mainestimate} are in stark contrast with the independent-variable case. For example, letting $(r_n)_{n\in \N}$ denote a sequence of i.i.d. Steinhaus random variables, one can show the asymptotically almost sure estimate 
\begin{equation}\label{independentest}
 \sup_{|t| \leq N^C} \Big|\frac1{\sqrt{N}} \sum_{n \leq N} r_n n^{it}\Big| \ll \sqrt{C \log N}
\end{equation}
for any fixed $C>0$ (see section \ref{contrast} below).
\end{enumerate}
\end{rmks}

\subsection{Background and related results}

A well-known (open) problem in analytic number theory is to determine the maximal size of the Riemann zeta function $\zeta$ in the critical strip, and on the critical line in particular. A conjecture of Farmer, Gonek, and Hughes \cite{FGH} asserts that
\begin{equation}\label{FGH_conj}
\max_{t\in[0,T]} |\zeta(\tfrac12 + it)| = \exp\Big( (\tfrac{1}{\sqrt{2}} + o(1)) \sqrt{\log T \log_2 T}\Big).
\end{equation}

Given the approximation
\[
\zeta(\tfrac12 + it) \sim \sum_{n \le T} \frac{1}{n^{1/2+it}}, \qquad t \in [T, 2T],
\]
one can view $\sum_{n \le T} X(n) n^{-1/2 + it}$ as a random model for $\zeta$ (or, more generally, for a Dirichlet $L$-function), which preserves the multiplicative nature of the summands. Aymone, Heap, and Zhao studied this model with $t = 0$, see in particular \cite[Corollary 1]{AHZ} compared with \eqref{FGH_conj}. For ease of exposition, we have chosen to work with the unweighted model \eqref{Ddef}.
\\

A related, but somewhat different problem, with close ties to random matrix theory, is to study the distribution of the maximum of $\zeta$ in \emph{short} intervals. Fyodorov, Hiary, and Keating \cite{fyodorov2012freezing} conjectured the correct scaling of the local maximum, as well as properties of the limiting distribution. In large part, the aforementioned conjecture was established in works of Najnudel \cite{Naj18}, Harper \cite{harper2019partition}, and Arguin, Belius, Bourgade, Radziwi\l{}\l{}, and Soundararajan \cite{Arg19, arguin2020fyodorov}. Arguin, Ouimet, and Radziwi\l{}\l{} \cite{Arg21} consider the the maximal size of the zeta function over short intervals of varying length.

Harper \cite{harper2013note} suggested the following random model to approximate the behaviour of $\log|\zeta|$ in a short interval on the critical line,
\[
W_T(h) = \sum_{p \le T} \frac{\Re{\big[X(p) p^{-ih}\big]}}{\sqrt{p}}, \qquad h \in I \subset \R.
\]
Here, the sum runs over primes, so that the $X(p)$ are independent Steinhaus random variables. One can further simplify the analysis by replacing the $X(p)$ with independent standard (complex) Gaussian random variables. For results about the maximum size of these random models, see Arguin, Belius, and Harper \cite{Arg17} (for intervals of fixed length), and Arguin, Dubach, and Hartung \cite{ArgDub21} (for intervals of varying length).




More generally, the distribution of partial sums of random multiplicative functions have been studied extensively. See, for instance, Basquin \cite{Bas}, Aymone, Fr\'{o}meta, and Misturini \cite{AFM} , Chatterjee and Soundararajan \cite{CS}, Harper, Nikeghbali, and Radziwi\l{}\l{} \cite{HNM}, Heap and Lindqvist \cite{HL}, Klurman, Shkredov, and Xu \cite{KSX}.

One might use a real-valued RMF to obtain a random counterpart to the Liouville function $\lambda$, or the Möbius function $\mu$ (if the RMF is non-zero just for squarefree values). In the non-random setting, the problem of establishing conditional and unconditional estimates for partial Möbius sums has  attracted the attention of numerous authors.  For instance,  starting with the work of Landau \cite{Lan}, various estimates of the form
$$M(x)=\sum_{n \leq x} \mu(n) \ll x^{1/2}  \exp((\log x)^{\theta} (\log \log x)^{\rho} )$$
have been shown to hold under the Riemann Hypothesis. Soundararajan \cite{Sound} established the above bound with exponents $\theta=1/2$ and $\rho=14$ (see also Balazard and De Roton \cite{Rotard}). In fact, one can prove even stronger estimates (including lower bounds) for $M(x)$ under the assumption of various far-reaching conjectures, see Ng~\cite{Ng}. Finally, in this context, we also mention the work of Maier and Sankaranarayanan \cite{MS}, which deals with more general, M\"{o}bius-like coefficients. \\ 

\subsection{Supremum over the real line}
It is natural to investigate the supremum of $D_N$ over the entire real line, that is to say
 $$\mathcal{M}_N=\E \sup_{t \in \R} |D_N(t)|,$$
in particular since this quantity determines the abscissa of uniform convergence $\sigma_u$ of the Dirichlet series $\mathcal{D}(s)=\sum_{n\geq 1} X(n)n^{-s}$. We recall that $\sigma_u$ is defined to be the infimum of those values $\sigma$ for which the series $\mathcal{D}(\sigma+it)$ converges uniformly over all $t \in \R$; it may be computed via the formula
$$\sigma_u=\limsup_{N \rightarrow \infty}  \frac{\log \mathcal{M}_N}{\log N}.$$
For independent variables $r_n$ it was shown that (see Lifshits and Weber \cite{LW} and Queff\'{e}lec \cite{Q})
$$\E \sup_{t \in \R} \Big| \sum_{n \leq N} r(n) n^{-\sigma+it}\Big| \asymp_{\sigma} \frac{N^{1-\sigma} }{\log N} $$
for $0 \leq \sigma \leq 1/2$ (here the implied constants in the upper and lower bound depend only on $\sigma$). However, in the case of the RMF $X(n)$ one obtains the trivial identity $\mathcal{M}_N=N$ (corresponding to $\sigma = 0$) from Bohr's correspondence. This simple fact will be explained in the final section of the paper.

\subsection{Notation}
The symbol $p$ is reserved for prime numbers and the expression $n \asymp x$ means $n \in [x/2,x]$. We will write $f \ll g$ or alternatively $f=O(g)$, if there exists an absolute constant $C$ such that $|f|\leq C |g|$. Often times we will add a subscript $f \ll_t g$ to emphasize the dependence of the implicit constant $C$ on the parameter $t$. The expression \textsl{natural parameter} refers to any quantity in $\N$.  The shorthand $\log_2(x)=\log \log x$ will be in use and we let $\omega(n)$ resp. $\Omega(n)$ denote the number of prime divisors of $n$, counted without resp. with multiplicity. Finally, the superscript $\flat$ will indicate a summation over squarefree variables while the symbol $\square$ indicates perfect square integers.

\subsection{Acknowledgements}
It is our pleasure to thank Ofir Gorodetsky, Adam Harper,  Herv\'{e} Queff\'{e}lec and Mikhail Sodin for helpful conversations and suggestions. In particular, Harper pointed out a nice approach, used in section 3, which simplifies an argument used in a previous version of this paper.
We thank the referee for a careful reading of this work, leading to an improved introduction, and sharper statement of the main result.
J.B. was supported by ERC Advanced Grant 692616. The research of A.N. was funded in part by ISF Grant 1903/18.

\section{Some divisor sums}

Let us begin with some notation and definitions. Recall the (usual) $k$-fold divisor function
$$
\tau_k(n) = \sum_{ \substack{ a_1, \dots, a_k \ge 1 \\ \: a_1 \cdots a_k = n }} 1.
$$
For $\alpha \in (0,1)$ we will call an integer $\alpha$-regular (of height $M$) if it belongs to the set
\begin{equation}\label{Gammadef} 
\Gamma_{\alpha}(M)=\left\{ m \leq M: \Omega(m) \leq (\log M)^{\alpha} \right\}.
\end{equation}
The complement of $\Gamma_{\alpha}(M)$, consisting of all $\alpha$-\textsl{irregular} numbers will be denoted by
\begin{equation}\label{Gammacdef}
\widetilde{\Gamma}_{\alpha}(M)= \left\{m \leq M: \Omega(m)> (\log M)^{\alpha}  \right\}.
\end{equation} 
Next we define two types of modified $k$-fold divisor functions. Given a real parameter $R \geq 2$, let us introduce 
$$\tau_{k,R}(n)=\sum_{ \substack{ a_1,\dots,a_k \leq R\\ a_1 \dots a_k= n }} 1, \qquad \qquad  \tau_{k,R;\alpha}(n)=\sum_{ \substack{ a_1,\dots,a_k \in \Gamma_{\alpha}(R)\\ a_1 \dots a_k= n }} 1.$$
 
\subsection{Upper bounds}
The goal of this section is to give a mean upper bound for the divisor function $\tau_k$ in moderately short intervals. Before stating the main result in section \ref{short}, we give some pointwise and mean-value estimates for divisor functions and binomial coefficients. Since we are working with large values of $k$ (relative to the length of the summation interval), some care is required.

\subsubsection{Preliminaries}
Let us first record the useful inequalities
\begin{equation}\label{optimal}
\sup_{x \geq 1} \ \left( \frac{a}{x} \right)^{x} \leq \exp(a/e), \qquad a>0
\end{equation}
and 
\begin{equation}\label{basictaubounds}
\tau_j(n)\tau_k(n)\leq \tau_{jk}(n), \qquad \tau_{j}(mn)\leq \tau_j(m) \tau_j(n).
\end{equation}
These last two inequalities hold for all $j,k \geq1$ and $m,n \geq 1$ and are easily verified, first at prime powers $n=p^e, m=(p')^v$ and then by extending multiplicatively. 

\begin{lem} For natural numbers $n>r \geq 1$ we have the bounds 
\begin{align}
 \left( \frac{n}{r} \right)^r & \leq { {n }\choose{ r} } \leq   \left( \frac{en}{r} \right)^{r} \label{binomial1}\\
		  \log { {n }\choose{ r} } & \leq n/3 , \qquad \qquad r/n \in[0.9,1]. \label{binomial2}
\end{align} 
\begin{proof}

The inequalities in \eqref{binomial1} are standard. To prove \eqref{binomial2}, we invoke the well-known estimate
$$ { {n }\choose{ r} } \leq \sqrt{\frac{n}{2 \pi r(n-r)}} \exp(n H(r/n)),$$
which is a straightforward consequence of Stirling's approximation. Here $H(x)=-x \log x -(1-x) \log(1-x)$ and a simple calculation reveals that $H(y) < 1/3$ in the range $y \in [0.9,1]$, yielding \eqref{binomial2}.
\end{proof}
\end{lem}

\begin{lem}\label{4divisor}
For any natural parameters $k,s \geq 1$,we have the uniform estimate  
\begin{equation}
\sum_{m \leq M} \tau_{k}(m)^s  \leq M (2 \log M)^{k^s-1}. \label{divisorpowersum} 
\end{equation}
Moreover, for all sufficiently large $k$ and any real $\sigma \in [9/5,2]$, we have that 
\begin{equation}
\sum_{m \geq 1} \frac{ \tau_{k}(m^2 )}{m^{\sigma} }  \ll k^{10 k^{2/\sigma}} \label{harmonicdivisorsquares}. 
\end{equation}
\begin{proof} For $s=1$, the proof of \eqref{divisorpowersum} can be found in (\cite[Lemma 3.1]{JAB}), where it was shown that
$$\sum_{m \leq M} \tau_{\ell}(m)  \leq M (2 \log M)^{\ell -1}.$$
To bound the general $s$-th power divisor sum, we apply \eqref{basictaubounds}  to find that $ \tau_{k}(n)^s \leq \tau_{k^s}(n)$ and the result follows.\\
Next we consider the weighted `square' sum \eqref{harmonicdivisorsquares}. The LHS  is given by the Euler product
$$\sum_{m \geq 1} \frac{ \tau_{k}(m^2)}{m^{\sigma}}=\prod_{p} \left(1+ \sum_{j \geq 1} \frac{\tau_k(p^{2j}) }{p^{\sigma j}} \right)=:\prod_{p} A_k(p). $$
We will give two estimates for $A_k(p)$, depending on whether $p < 10k$ or not. From \eqref{binomial1} and \eqref{binomial2} we first find that
\begin{align}\label{akp}
A_k(p)&=1+\sum_{j \geq 1} p^{-\sigma j} { {2j+k-1 }\choose{ 2j} } \\
& \leq 1+ \sum_{j \leq 5 k}  \left( \frac{e(2j+k)}{2j \cdot p^{\sigma/2}}\right)^{2j} + \sum_{j \geq 5 k }p^{- \sigma j} \exp \left( \frac{2j+k}{3} \right). \notag
\end{align}
Since $\max_{\ell \geq1} (11 e k/ \ell p^{\sigma/2 } )^{\ell} \leq \exp(11 k/p^{\sigma/2})$ and $(2j+k)/3\leq 2j/2$ in the range $j \geq 5k$, we get that
\begin{align*}
A_k(p) & \leq 1+ \sum_{j \leq 5k}  \left( \frac{11ek}{2j \cdot p^{\sigma/2} }\right)^{2j} + \sum_{j \geq 5k}  \left( \frac{\sqrt{e}}{p^{\sigma/2}  }\right)^{2j} \\
& \leq 7+ 5 k \exp \left(\frac{11k}{  p^{\sigma/2}  } \right)
\leq 12 k \exp \left(\frac{11k}{  p^{\sigma/2}  } \right).
\end{align*}
In the last line we used the bound $\sqrt{e}/ p^{\sigma/2} \leq 9/10$. The above estimate will be useful when $p^{\sigma/2} < 7k$. On the other hand, when $p^{\sigma/2} \geq 7k$, we use the inequality $(2j+k)/2j \leq 2k$ (which holds for all $k,j \geq 1$) and proceed as in the first line of \eqref{akp} to get that
\begin{align*}
A_k(p) & \leq 1+ \sum_{j \geq 1}  \left( \frac{e(2j+k)}{2j \cdot p^{\sigma/2}}\right)^{2j} \leq 1+ \sum_{j \geq 1}  \left( \frac{2ek}{p^{\sigma/2}}\right)^{2j}  \\
& \leq 1+ \frac{(6k/ p^{\sigma/2})^2}{1-(6k/ p^{\sigma/2})^2}\leq 1+4\left( \frac{6k}{ p^{\sigma/2}} \right)^2 .
\end{align*}
Combining the two estimates for $A_k(p)$ and Chebychev's upper bound for the density of primes, we now gather that
\begin{align*}
\prod_p A_k(p) & \leq \prod_{p^{\sigma/2} \leq 7 k} \left[ 12 k \exp(11k/p^{\sigma/2} ) \right]
 \prod_{p^{\sigma/2} \geq 7 k }  \bigg[ 1+4 \Big( \frac{6k}{p^{\sigma/2}}  \Big)^2 \bigg] \\
& \leq (12 k)^{(7k)^{2/\sigma}} \exp \bigg( 11k \sum_{p^{\sigma/2} < 7k} \frac{1}{p^{\sigma/2}} + B_1 \sum_{p^{\sigma/2}\geq 10k } \frac{k^2}{p^{\sigma} } \bigg)\\
& \leq (12 k)^{ 9k^{2/\sigma}}  \exp \left( B_2  \ k^{2/\sigma} \log_2 k+ B_3 \  k^{2/\sigma}\right)  \ll k^{10 k^{2/\sigma }}
\end{align*}
for $k>k_0$ sufficiently large. The precise values of the absolute constants $B_j>0$, appearing in the last two lines, are unimportant.
\end{proof}
\end{lem}

\subsubsection{Divisor sums in short intervals}\label{short}
The key result of this section, stated in Proposition~\ref{mainuppershort} below, deals with short divisor sums of the form $\sum_{n \in [X,X+Y]} \tau_k(n) $. In view of our specific applications it is important that we let $k$ grow faster than $(\log X)^{1-\varepsilon}$. For such large values of $k$ it will be convenient to work with $(1-\varepsilon)$-regular integers $n$, that is to say $n \in \Gamma_{1-\varepsilon}(X)$.

Before moving on to the proposition, we require two more ingredients. First we will need the classical bound of Hardy-Ramanujan which controls the number of integers with an unusually large amount of prime divisors. In its original form \cite{HR}, the theorem asserts the existence of a constant $c>0$ such that
\begin{equation}\label{HRnu}
\pi_{\nu}(x):=\left|\left\{ n \leq x: \omega(n)=\nu \right\}\right|\ll \frac{x}{\log x} \frac{(\log_2 x+c)^{\nu -1}}{(\nu -1)!},  
\end{equation}
for any natural number $\nu$ and any $x\geq 3$. 
The second ingredient is the following squarefree version of the main proposition.

\begin{lem}\label{squarefreemainprop}
For any exponent $\sigma \in [1/2,1]$ and any pair of parameters  $X \geq 10$ and $X^{\sigma} \leq Y \leq X$,  we have the estimate 
\begin{equation}\label{shortsfreediv}
 \sum_{n \in [X,X+Y]  }^{\flat} \tau_k(n) \ll Y (\log X)^4 \exp(2 k^{1/\sigma} \log_2 X),
\end{equation}
provided that $k \leq \log X$.
\begin{proof} 
Let us first give a short-interval version of \eqref{HRnu}. Given any squarefree $n\in [X,X+Y]$ with $\nu$ prime factors, we let $d| n$ be the divisor formed by the  product of the $\nu_{\sigma}=\lfloor \sigma \nu \rfloor $ smallest prime divisors of $n$; clearly $d \leq 2X^{\sigma}$.  As a consequence of \eqref{HRnu} we find that 
\begin{align}\label{shortHRnu}
 \sum_{\substack{ n \in [X,X+Y]\\ \omega(n) = \nu }} 1 & \leq \sum_{ \substack{ d \leq 2X^{\sigma}\\ \omega(d)=\nu_{\sigma}    } }\  \sum_{r \in [ \frac{X}{d },\frac{X+Y }{d } ]} 1 \notag \\
& \ll Y \sum_{ \substack{ d \leq 2X^{\sigma}\\ \omega(d)=\nu_{\sigma}    } } \frac{1}{d } \ll Y \frac{(\log_2 X +c)^{\nu_{\sigma} } }{(\nu_{\sigma} -1)!}.
\end{align}
The bound in the last line follows after a simple dyadic decomposition of the interval $[1,X^{\sigma}]$.\\ 
Proceeding with the treatment of \eqref{shortsfreediv},  we apply \eqref{shortHRnu}, keeping in mind the pointwise bound $\tau_k(n)\leq k^{ \omega(n)}$ (valid for squarefree $n$) and find that the
 LHS of \eqref{shortsfreediv} is no greater than
\begin{align*}
 \sum_{\nu \leq \log X }\sum_{ \substack {    n \in [X,X+Y]   \\ \omega(n)=\nu  }}^{\flat} k^{\nu } &\ll (\log X) Y \max_{ \nu \leq \log X  }   \frac{ k^{\nu }  (2 \log_2 X)^{\nu_{\sigma} }}{ (\nu_{\sigma} -1)!} \\
& \leq (\log X)^4 \, Y \max_{ \nu_{\sigma}  \geq 1  }   \frac{   (2k^{1/\sigma} \log_2 X)^{\nu_{\sigma}   }}{ \nu_{\sigma} !}.
\end{align*}
Inserting the lower bound $\nu_{\sigma}  ! \geq (\nu_{\sigma} /e)^{\nu_{\sigma} }$ into the last line and applying \eqref{optimal}, we easily retrieve \eqref{shortsfreediv}.

\end{proof}
\end{lem}

\begin{rmk}\label{sigmaremark}
It is important to point out the unnatural expression $k^{1/\sigma}$ appearing on the RHS of \eqref{shortsfreediv}. In particular, when $k$ is large and $[X,X+Y]$ is a short interval (e.g. $\sigma=1/2$), the estimate is very poor. This loss of accuracy stems from the  bound \eqref{shortHRnu} and it would be interesting to determine whether the factor $1/(\nu_{\sigma}-1)!$ may be sharpened to $1/(\nu -1)!$ in the short-interval setting. Fortunately, for our applications, it will be enough to apply Lemma \ref{squarefreemainprop} with values of $\sigma$ approaching $1$. 
\end{rmk}

\begin{prop}[Divisor sums in short intervals]\label{mainuppershort}
Let $\alpha \in (0,1)$ and $\sigma \in [9/10,1]$ be fixed. Then for any sufficiently large $M \geq M_{\alpha} \geq 1 $, any integer $1 \leq k \leq \log M$ and parameter $M^{\sigma} \leq H \leq M$, we have the estimate
\begin{equation}\label{shortdivisor}
\sum_{\substack{ m \in [M,M+H] \\ m \in \Gamma_{\alpha}(M) }  } \tau_{k}(m)  \ll_{ \alpha}  H   \exp(20 k^{1/\sigma^2} \log_2 M). 
\end{equation}
\begin{proof} Let us first assume that $\sigma \leq \sigma_0:= 1-1/100 \log_2 M$.  We observe that each natural number $m$ admits a factorisation $m=\widetilde{m}\cdot m'$ where $\widetilde{m}$ is its largest square divisor and $m'$ is squarefree.  In view of this factorisation and the pointwise bound $\tau_k(m) \leq k^{\Omega(m)} \le H^{(1-\sigma)/2} $, which holds for any $m \in \Gamma_{\alpha}(M)$, and $M \ge M_\alpha$, we may separate the sum on the LHS of \eqref{shortdivisor} to get that   
\begin{align*}
\sum_{\substack{ m \in [M,M+H] \\ m \in \Gamma_{\alpha}(M) }  }  \tau_{k}(m) \ll_{\alpha} & 
\sum_{\substack{ m=\widetilde{m} m' \in [M,M+H ]\\H^{1-\sigma} \leq \widetilde{m} \leq M } }  H^{\frac{1-\sigma}{2} } 
+ \sum_{ \substack{ \widetilde{m} < H^{1-\sigma} \\  m \in \Gamma_{\alpha}(M)} } \tau_k(\widetilde{m}) 
\sum_{m' \in [\frac{M}{\widetilde{m} }, \frac{M+H}{\widetilde{m}  }  ]}^{\flat} \tau_k(m')\\
				&=: H^{\frac{1-\sigma}{2} } \  \mathcal{T}_1+\mathcal{T}_2.
\end{align*}
Let us first consider $\mathcal{T}_1$.  Since $H\geq M^{\sigma}$, we have that
\begin{align*}
\mathcal{T}_1 & \leq    \sum_{\substack{ m=\widetilde{m} m' \in [M,M+H ]\\H^{1-\sigma} \leq \widetilde{m} \leq H} } 1
+ \sum_{\substack{ m=\widetilde{m} m' \leq 2M \\H \leq \widetilde{m} \leq M} } 1  \\
&  \ll  \sum_{  \substack{ H^{1-\sigma} \leq \widetilde{m} \leq H  }}^{\square}  \frac{H}{\widetilde{m}} + 
\sum_{ \substack{ \widetilde{m} \geq H  }}^{\square}  \frac{M}{\widetilde{m} }    \ll  
 \frac{H}{H^{\frac{1-\sigma}{2} }  }+ \frac{M}{H^{1/2} }    \ll \frac{H}{H^{\frac{1-\sigma}{2} }  }.
\end{align*}
Here, the superscript $\square$ indicates a summation over perfect squares. \\
To treat  $\mathcal{T}_2$, we observe that $H/ \widetilde{m}  \geq (M/ \widetilde{m} )^{\sigma^2}$ whenever $\widetilde{m} \leq H^{1-\sigma}$.  Combining \eqref{harmonicdivisorsquares} and \eqref{shortsfreediv} it follows that
$$\mathcal{T}_2 \ll H  (\log M)^4 \exp(2k^{1/\sigma^2} \log_2 M)  \sum_{  \widetilde{m} < H^{1-\sigma} }^{\square} 
\frac{ \tau_k(\widetilde{m}) }{\widetilde{m}} \ll H  \exp(16 k^{1/\sigma^2} \log_2 M).  $$
We may now collect the estimates for $\mathcal{T}_1$ and $\mathcal{T}_2$, concluding that the LHS of \eqref{shortdivisor} is $O
_{\alpha}(H  \exp(16 k^{1/\sigma^2} \log_2 M))$, provided that $\sigma \leq \sigma_0=1-1/100 \log_2 M$. \\
To conclude the argument, all that remains is to consider values $\sigma \in [\sigma_0,1]$. In this case, it is enough to split the range $[M,M+H] $ into shorter intervals, each of length $M^{\sigma_0}$,  and apply the estimate proven just above, together with the straightforward inequality $k^{1/\sigma_0^2} \leq \frac{5}{4} k \leq \frac{5}{4} k^{1/\sigma^2}$.
\end{proof}

\end{prop}

The final estimate in our series of  mean upper bounds for  divisor functions concerns $\alpha$-irregular numbers.

\begin{lem}[Irregular divisor sums]\label{irregest} Given natural parameters $k,M \geq 10$ and any $\alpha \in (0,1)$, we have that  
\begin{equation}\label{irregdivbound}
\sum_{m \in \widetilde{\Gamma}_{\alpha}(M)} \tau_{k}(m)  \ll_{\alpha} M \exp(-(\log M)^{\alpha}/ 4 \log_2 M), 
\end{equation}
provided that $k \leq (\log M)^{\alpha}/ (\log_2 M)^3$.
\begin{proof}
As in the proof of Proposition \ref{mainuppershort}, we may write $m=\widetilde{m}\cdot m'$ with $m'$ squarefree. Using this factorisation, we choose a threshold parameter $Y=\exp((\log M)^{\alpha} /2)$ and consider the cases $\widetilde{m} \geq Y$ and $\widetilde{m} < Y$ separately. Observe that for any $m=\widetilde{m}\cdot m' \in \widetilde{\Gamma}_{\alpha}(M)$ satisfying $\widetilde{m} < Y$, we necessarily have that $\omega(m') \geq (\log M)^{\alpha}/4$. We gather that
\begin{align*}
\sum_{m \in \widetilde{\Gamma}_{\alpha}(M)} \tau_{k}(m)& \leq   \sum_{Y \leq  \widetilde{m} \leq M}^{\square} \tau_k(\widetilde{m}) \sum_{ m' \leq M/  \widetilde{m} } \tau_{k}(m') 
+ \sum_{ \widetilde{m} < Y}^{\square} \tau_k(\widetilde{m}) \sum_{ \substack {m' \leq M/  \widetilde{m} \\ \omega(m') \geq (\log M)^{\alpha}/4 }} \tau_{k}(m') \\
& =:  \mathcal{K}_1+  \mathcal{K}_2.
 \end{align*}
To bound $ \mathcal{K}_1$, we apply \eqref{harmonicdivisorsquares},  taking $\sigma=2-(\log_2 M)^{-1}$. Using \eqref{divisorpowersum} to deal with the inner-most sum,  we find that 
\begin{align*}
 \mathcal{K}_1 & \leq M (2 \log M)^{k} \sum_{ Y \leq  \widetilde{m} \leq M}^{\square} \frac {\tau_k(\widetilde{m}) }{\widetilde{m} } 
 \leq M (2 \log M)^{k} Y^{-\frac{2-\sigma}{2}} \sum_{ Y^{1/2} \leq  r\leq M^{1/2} } \frac {\tau_k(r^2) }{r^\sigma } \\
& \ll M(2 \log M)^{k}  \exp \left( -\frac{(\log M)^{\alpha}}{ 2 \log_2 M} \right) k^{10 k^{ 2/ \sigma} } \ll M \exp(-(\log M)^{\alpha}/ 3 \log_2 M) .
\end{align*}
To treat $\mathcal{K}_2$ we apply \eqref{harmonicdivisorsquares}  once again,  this time setting $\sigma=2$. In view of \eqref{HRnu} we gather that
\begin{equation*}
 \mathcal{K}_2 \leq M \sum_{ \widetilde{m} \leq M}^{\square} \frac {\tau_k(\widetilde{m}) }{\widetilde{m} } \sum_{(\log M)^{\alpha}/4 \leq \nu \leq \log M} \left( \frac{3 k \log_2 M}{\nu} \right)^{\nu}
 \leq M \log M \exp(-(\log M)^{\alpha}).
\end{equation*}
Collecting the estimates for $ \mathcal{K}_1$ and $ \mathcal{K}_2$, we get \eqref{irregdivbound}.
\end{proof}
\end{lem}

\subsection{A lower bound for the second moment of $\tau_{k,R;\alpha}$ }
In order to furnish a lower bound for the sum $\sum_{n \leq N} \tau_{k,R;\alpha}(n)^2$, we will restrict the values of $n$ to a suitable subset of integers. Given any natural $1 \leq \nu \leq \log N/ (\log_2 N)^2$, set
\begin{align*}
L&=(\log N)^3, \qquad  L'=\lceil L/ 3 \log L \rceil, \\
Y&=N/L^{\nu}, \qquad  Y'=\lceil Y/ 3 \log Y \rceil. 
\end{align*}
We may then define the collections 
\begin{equation*}
\mathcal{P}(L; \nu)=\left\{q=p_1\cdots p_{\nu} : p_j \asymp L \ \text{ for all $j$, with } p_j \text{ distinct} \right\}
\end{equation*}
and
\begin{equation*}
\mathcal{G}(N; \nu) =\left\{q\cdot p': q \in \mathcal{P}(L; \nu),\  p' \asymp Y \right\}.
\end{equation*}
Finally, let  
\begin{equation*}
\mathcal{A}(N;k, \nu) =\left\{A=n_1\cdots n_k: \forall i \neq j \  \gcd(n_i,n_j)=1, \text { and }  n_j \in \mathcal{G}(N; \nu) \right\}.
\end{equation*}
We also record the following weak version of Stirling's approximation: for any natural $r \geq 1$ we have that  
$$\sum_{n \leq r} \log n = \log (r!) \geq r \log r -r. $$

\begin{lem}\label{mainlower} Let $\alpha \in (0,1)$ be given, suppose that $N \geq N_0$ is sufficiently large and $(\log_2 N)^3 \leq k \leq (\log N)^{\alpha}$. We have the lower bound
\begin{equation}\label{lowerdivisork}
\sum_{A \leq N^k} \tau_{k,N;\alpha} (A)^2 \gg N^k\exp \left( \frac{k^2}{200 \log_2 N} \right).
\end{equation}
\begin{proof}
Let $\nu<k$ be a large parameter, to be chosen later. It is enough to restrict the LHS of \eqref{lowerdivisork} to values $A \in \mathcal{A}(N;k, \nu) $ and give a lower bound for the resulting divisor sum. To this end we first estimate the cardinality of $ \mathcal{A}(N;k, \nu) $ from below. Since each element $A$ in the set is obtained by choosing $\nu k$ distinct primes in the interval $[L/2,L]$ and $k$ distinct primes in the interval $[Y/2,Y]$, we gather that 
\begin{align*}
|\mathcal{A}(N;k, \nu)| &\geq \binom{L'}{\nu k} \binom{Y'}{k} 
	 \geq \frac{ (L^{\nu} Y)^{k} }{(3 k \log Y)^k (3 \nu k \log L)^{\nu k}} \\
	& \geq \frac{ N^{k} }{(3 k \log N)^k (3 \nu k \log L)^{\nu k}}.
\end{align*}
Moreover, for each $A \in \mathcal{A}(N;k, \nu)$, the number of ways to obtain a factorisation $A=n_1\cdots n_k$ with (pairwise coprime) $n_j \leq N$ is at least
\begin{align*}
  \prod_{j=0}^{k-1} { {\nu (k-j) }\choose{ \nu} } \geq  \exp \left( \nu \sum_{j=0}^{k-1} \log(k-j)   \right) \geq  \left( \frac{k}{e} \right)^{k \nu}.
\end{align*}
The two previous estimates combined, yield
\begin{align*} 
\sum_{A \leq N^k} \tau_{k,N;\alpha} (A)^2 & \geq \sum_{A \in \mathcal{A}(N;k, \nu) }  \left[  \Big( \frac{k}{e} \Big)^{k \nu} \right]^2\\
& \gg \frac{ N^k }{3 (k \log N)^k}\left( \frac{k^2}{3 e^2 \nu k \log L} \right)^{k \nu}.
\end{align*}

Inserting the choice of parameter $\nu=\lfloor k/30 \log L\rfloor $ into the lower bound just above, we easily retrieve \eqref{lowerdivisork}.
\end{proof} 
\end{lem}

\section{Proof of the main result}

As in \cite{JAB}, we will use a moment method to control the size of $\sup_{|t| \leq N^{C(N)}} |D_N(t)| $; often times we will use the notation $T=N^{C(N)}$. Although in many arguments of this section, the quantity $C>0$ will be allowed to grow/decay with $N$ in an arbitrary fashion, there are crucial estimates (such as \eqref{upperErem} in Lemma \ref{supupperlemma} below) which require the additional assumption that $C=C(N)$ satisfy \eqref{Cgrowth} as $N \rightarrow \infty$.

\begin{definition}
For any real parameter $T\geq 1$, we define the random variable given by the $2k$-th moment
\begin{equation}\label{Mkdef}
M_k=M_k(T)= \int_{-T}^{T} |D_N(t)|^{2k} \ dt.
\end{equation}
\end{definition}
Our starting point for both the upper and lower bound in \eqref{mainestimate} is the evaluation of $\E[M_k]$. Using the convenient notation ${\bf n}=(n_1,...,n_k)$ for $k$-tuples of integers, together with the identity
\begin{equation*}
\E[X(p)^{r_1} \overline{X(p)^{r_2}}]=\delta_{r_1,r_2} 
\end{equation*}
for any natural powers $r_1, r_2$, we gather that 
\begin{align}\label{expectationmk}
\E[M_k]& =N^{-k}\int_{-T}^{T} \E \Big[ \sum_{{\bf n},{\bf m} \in [1,N]^k } X(n_1 \cdots n_k) \overline{X(m_1 \cdots m_k)} \Big]
\left( \frac{n_1 \cdots n_{k}}{m_1 \cdots m_{k}} \right)^{it } \ dt   \notag \\
&=N^{-k}\int_{-T}^{T} \sum_{\substack{{\bf n},{\bf m} \in [1,N]^k \\ n_1 \cdots n_{k}=m_1 \cdots m_{k}  }} \left( \frac{n_1 \cdots n_{k}}{m_1 \cdots m_{k}} \right)^{it } \ dt \notag \\
&= 2TN^{-k} \sum_{A \leq N^k} \tau_{N,k}^2 (A).
\end{align}

\subsection{ The upper bound}\label{sectionupper}
\subsubsection{Moment estimates}
In order to estimate $ \sup_{|t| \leq N^C} |D_N(t)|$ from above, it will be enough to bound $|D_N(t)|$ pointwise, provided that we can do so with high probability. We are grateful to Adam Harper for suggesting this approach since it greatly simplifies the argument we gave in a previous version of this paper.\\ 
First we give some notation. Let $\mathscr{I}= \left\{I_j\right\}_{j \leq J}=\left\{ [a_j,b_j )\right\}_{j \leq J}$ be a collection of intervals such that $J=O(T )$ and
\begin{equation}\label{Iconditions}
[-T,T)= \dot{\bigcup}_{j \leq J} I_j, \qquad \qquad |I_j| \asymp 1, \qquad \text{for all } j=1,\dots,J.
\end{equation}
Given any Dirichlet polynomial $d(t)$, and natural parameters $r,\ell \geq 0$, let us write
\begin{equation}\label{Sobmomentdef}
\mathcal{D}_{\ell}^{(r)}(I_j)=\int_{I_j} |d^{(r)}(t)|^{2\ell} \ dt.
\end{equation}
Assuming that $|d(t)|$ takes its maximum at $t_ j \in I_j$ for each interval $I_j \in \mathscr{I}$, we have that  
\begin{equation}\label{meantosup}
\sup_{t \in I_j} |d(t)|= \bigg| d(t_j) + \int_{a_j}^{t_j} d'(t) \ dt \bigg| \leq |d(a_j)|+   \int_{I_j} |d'(t)| \ dt .
\end{equation}

In certain settings the following crude alternative will be of use (cf. \cite[Lemma~4.2.]{JAB}).
\begin{lem}
For any $k \geq 1$, $T \geq 1$ and Dirichlet polynomial  $d(t)=\sum_{n \leq N} a_n n^{it}$ of length $N \geq 3$, we have that
\begin{equation}\label{sup}
 \sup_{|t|\leq T}|d(t)|  \ll \left( \log N  \|a\|_1 \ \mathcal{D}_k^{(0)}([-2T,2T])\right)^{\frac{1}{2k+1}} .
\end{equation}
\begin{proof}
Let $H = \sup_{|t|\leq T} |d(t)| $ and define $S=\|a\|_1 \log N \geq  \sup_{ |t| \leq 2T } |d'(t)|$.  Thus if $|d(t)|$ achieves its maximum at $t = t_0\in [-T,T]$,  we gather that $|d(t_0 + t)| \geq H - S|t|$, and hence $|d(t_0 + t)| \geq H/2$ whenever $|t| \leq \frac{H}{ 2  S}$. 
Since $H/2S \leq 1 \leq T$, we find that 
$$\frac{H}{S}(H/2)^{2k}\leq \int_{  [t_0-H/2S, t_0 +H/2S ]  }  |d(t)|^{2k}\ dt  \leq \int_{-2T}^{2T} |d(t)|^{2k}\,dt.$$
\end{proof}
\end{lem}

The purpose of the next lemma is to furnish an upper bound for the expectation of $ \sup_{t \in I_j} |D_N(t)|$. 
We also record a variant of the estimate for the `remainder' polynomial $\widetilde{D}_N^{\alpha}(t)$ which is defined as follows. Let $\alpha \in (0,1)$ and set 
\begin{equation}\label{DNalpha}
\widetilde{D}_N^{\alpha}(t)= \frac1{\sqrt{N}} \sum_{\substack{ n \leq N \\ n \in \widetilde{\Gamma}_{\alpha}(N)}} X(n) n^{it},
\end{equation}
where $\widetilde{\Gamma}_{\alpha}(N)$ is the complement of the set $\Gamma_{\alpha}(N)$, as given in \eqref{Gammacdef}. 

\begin{lem}\label{supupperlemma}
a) Let $N \geq 100$ and $1 \leq k \leq \frac{1}{2} \log N$ be given. Then for any $C>0$ and any partition $\mathscr{I}$ satisfying \eqref{Iconditions}, the estimate
\begin{equation}\label{upperEregCsmall}
\E \bigg[\sup_{ t \in I_j} |D_N(t)|^{2k} \bigg] \ll   4^k  ( \log N)^{2(k+1)^2 }
\end{equation}
holds for each $j \leq J$.\\
 b) Let $\gamma \in (0,1)$ be the upper exponent in \eqref{Cgrowth} and suppose that $\alpha \in (\frac{1}{2}(1+ \gamma),1)$, $N \geq N_{\alpha}$ is sufficiently large and $k \asymp 20 (\log N)^{1+\gamma-\alpha} \log_2 N$. Assuming that $C=C(N)$ satisfies \eqref{Cgrowth}, we have that
\begin{equation}\label{upperErem}
\E \bigg[\sup_{|t| \leq N^C} |\widetilde{D}_N^{\alpha}(t)| \bigg] \ll_{ \alpha, \gamma} 1.
\end{equation}
\begin{proof}

$\mathit{a)}$  To establish \eqref{upperEregCsmall}, we will first give the necessary estimates for an application of \eqref{meantosup}.
To begin with, we need to treat the moments of the derivative 
\[
D_N'(t) = \frac{i}{\sqrt{N}} \sum_{n \leq N}   X(n)   n^{it} \log n.
\]
Let us write $T=N^C$. Recalling the notation \eqref{Sobmomentdef} (with $d(t)=D_N(t)$), it follows from  \eqref{divisorpowersum} that for any $\ell \geq 1$ and $j \leq J$
\begin{align}
 \E [ \mathcal{D}_{\ell}^{(1)}(I_j) ] &  =N^{-\ell} \int_{I_j}  \sum_{\substack{{\bf n},{\bf m} \in [1,N]^{\ell} \\ n_1 \cdots n_{\ell}=m_1 \cdots m_{\ell}  }} \prod_{j \leq \ell} \left( \log n_j \log m_j  \right) d t \notag\\
& \leq |I_j| N^{-\ell} (\log N)^{2 \ell } \sum_{A \leq N^{\ell} } \tau_{\ell}(A)^2 \leq (2 \log N^{\ell})^{(\ell+1)^2} \label{upperEresult}
\end{align}
and the exact same argument yields (recall the notation $I_j=[a_j,b_j)$)
\[ \E[|D_N(a_j)|^{2 \ell}] \leq    (2 \log N^{\ell})^{\ell^2}.  \]
We now let $S_j= \sup_{t \in I_j} |D_N(t)| $ and invoke \eqref{meantosup}. Given any $1 \leq k \leq \frac{1}{2} \log N$ we may apply H\"{o}lder's inequality, together with the above estimates, to find that
\begin{equation}
\E \left[S_j^{2k} \right]    \leq 4^k (\E [ \mathcal{D}_{k}^{(1)}(I_j) ] +  \E[|D_N(a_j)|^{2 k}] )
 \ll   4^k  ( \log N)^{2(k+1)^2 },   \label{Holderrep}
\end{equation} 
which recovers \eqref{upperEregCsmall}. \\
$\mathit{b)}$ The corresponding estimate for $\widetilde{D}_N^{\alpha}(t)$,  that is to say \eqref{upperErem}, is obtained in the same way as \eqref{upperEresult},  the only difference being that the summation variable $n$ runs over the set $\widetilde{\Gamma}_{\alpha}(N)$ which is very sparse. To be precise, we will assume that $k \asymp 20 (\log N)^{1+\gamma-\alpha}\log_2 N$ and then define $\beta \in (0,1)$ implicitly by way of the identity $k (\log N)^{\alpha}=(k \log N)^{\beta}$. As a consequence we have that $\Omega(A) \geq k (\log N)^{\alpha} \geq (\log (N^k))^{\beta}$
for any integer $A \in \widetilde{\Gamma}_{\alpha}(N)^k$. We now set $T=N^{C(N)}$ and proceed with a direct computation of the $2k$-th moment, together with an application of \eqref{basictaubounds} and Lemma \ref{irregest}. Recalling the notation \eqref{Sobmomentdef} once again (this time $d(t)=\widetilde{D}_N^{\alpha}(t)$), we find that
\begin{align}\label{CS}
\E [\mathcal{D}_{k}^{(0)}([-2T,2T])]
& \leq 4TN^{-k} \sum_{A \in\widetilde{\Gamma}_{\alpha}(N)^k  } \tau_k(A)^2  \notag \\
& \leq 4T N^{-k} \sum_{A \in\widetilde{\Gamma}_{\beta}(N^k)  } \tau_{k^2}(A) \notag \\
& \ll_{\alpha} T \exp(- k(\log N)^{\alpha}/8 \log_2 N ). 
\end{align} 
It should be noted that Lemma \ref{irregest} is applicable thanks to the inequality $1+\gamma-\alpha< \alpha$ which implies that $k^2 < k(\log N)^{\alpha}/  (\log_2 (N^k))^3= (\log (N^{k}))^{\beta}/  (\log_2 (N^k))^3$. Combining \eqref{CS}, \eqref{sup} and H\"{o}lder's inequality with the fact that $NT $ is dwarfed by $\exp( k( \log N)^{\alpha}/8 \log_2 N)$, we find that
\begin{align*}
\E \bigg[\sup_{|t| \leq N^C} |\widetilde{D}_N^{\alpha}(t)| \bigg] & \leq (N^{1/2} \log N)^{\frac{1}{2k+1} } \E [\mathcal{D}_{k}^{(0)}([-2T,2T])]^{\frac{1}{2k+1}  } \\
 &\leq  ( NT \exp(- k(\log N)^{\alpha}/8 \log_2 N ) )^{\frac{1}{2k+1}}  \ll_{ \alpha, \gamma} 1,
\end{align*}
as desired. 
\end{proof}
\end{lem}

\subsubsection{ Concluding the proof of the upper bound in Theorem \ref{mainprop}} 
Let us assume that $C(x)$ satisfies the estimates
$$\frac{(\log_2 x)^2}{\log x} \leq  C(x) \leq \log x$$
for all large $x$ and suppose that $N$ is a sufficiently large natural number. Now let $\left\{I_j\right\}_{j \leq J}$ be any partition satisfying \eqref{Iconditions}, fix any $\varepsilon >0$ and and set 
\[  k=\lfloor (C(N) \log N / \log_2 N)^{1/2} \rfloor -1, \qquad \qquad \lambda=\left (\tfrac{3}{4} +\varepsilon \right)(C(N) \log N  \log_2 N)^{1/2}. \]
Applying \eqref{upperEregCsmall}, we find that
\begin{equation*}
\E \bigg[ \sup_{t \in I_j} |D_N(t)|^{2k} \bigg] \leq 4^k \exp \left(  2 C(N) \log N  \right)
\end{equation*}
and,  as a result, we gain sufficiently strong control of the unlikely events 
$$E_{j}:  \left\{ \sup_{t \in I_j} |D_N(t)| \geq \exp(2 \lambda) \right\}. $$ 
Indeed, a straightforward application of Chebychev's inequality reveals that $\mathbb{P}(E_j)=C_{\varepsilon} \ o(1/T)$ for some $C_{\varepsilon}>0$ depending only on $\varepsilon>0$ and all that remains is to take the union bound $\mathbb{P}(\cup_{j \leq J} E_j)=C_{\varepsilon}o(1),$ recovering the upper bound in \eqref{mainestimate}.  It should also be noted that the RHS of \eqref{mainestimate} exceeds the trivial bound $\mathcal{S}(N,C) \leq N$ when $C(N)\geq \log N$.

\subsection{ The lower bound}\label{sectionlower}
In order to establish the lower bound in \eqref{mainestimate}, we would like to show that the $2k$-th moment $M_k$ concentrates around its mean by controlling the variance.  However, to avoid technical difficulties we will need to work with the following setup.  Let $\gamma \in (0,1)$, fix a value $\alpha \in (\frac{1}{2}(1+ \gamma),1)$ and let us first remove from $D_N(t)$, the remainder $\widetilde{D}_N^{\alpha}(t)$ defined in \eqref{DNalpha}.  Thanks to the estimate \eqref{upperErem}, we know that asymptotically almost surely
\begin{equation}\label{badsmall}
\sup_{|t| \leq N^{C(N)} } |\widetilde{D}_N^{\alpha}(t)| \ll_{\gamma, \alpha} \log N
\end{equation}
as $N \rightarrow \infty$. As a result, it will be enough to deliver an almost sure lower bound for the supremum of the `main part' 
$$D_N^{\alpha}(t)=D_N(t)-\widetilde{D}_N^{\alpha}(t)= \frac1{\sqrt{N}} \sum_{\substack{ n \leq N \\ n \in \Gamma_{\alpha}(N)}} X(n) n^{it}.$$
We recall the set of $\alpha$-regular integers $\Gamma_{\alpha}(M)=\left\{ m \leq M: \Omega(m) \leq (\log M)^{\alpha} \right\}$  appearing in the last line and, accordingly, consider the modified moments
\begin{equation}\label{alphamoment}
M_{k,\alpha}=M_{k,\alpha}(T)=\int_{-T}^{T} |D_N^{\alpha}(t)|^{2k} \ dt.
\end{equation} 
The same calculation as \eqref{expectationmk} gives the evaluation
\begin{equation}\label{expectationmkalpha}
\E[M_{k,\alpha}]= 2TN^{-k} \sum_{A \leq N^k} \tau_{k,N;\alpha}^2 (A)
\end{equation} 
and thus, in view of Lemma \ref{mainlower}, we are left with the task of bounding $\Var[M_{k,\alpha}]$.

\begin{prop} Let $\alpha \in (0,1)$ be given. Then for any $T \geq 1$,  all sufficiently large $N \geq 1$, and any $k \leq \log N$,  we have the estimate
\begin{equation}
\Var[M_{k,\alpha}] \ll T^{2-\rho}    \exp(200k^2 \log_2 N). 
\end{equation}
Here we have used the notation $\rho=(1000 \log_2 N)^{-1}$. 
\begin{proof}
To get a handle on the variance we first write 
\begin{align*}
\E[M_{k,\alpha}^2]& =N^{-2k}\int_{-T}^{T} \int_{-T}^{T}  \E \bigg[ \sum_{{\bf n},{\bf m} \in [1,N]^k }^{\sharp} \sum_{{\bf n'},{\bf m'} \in [1,N]^k }^{\sharp} X(n_1 \cdots n_k) \overline{X(m_1 \cdots m_k)} 
 \\
& \times  X(n_1' \cdots n_k') \overline{X(m_1' \cdots m_k')} \bigg]
\left( \frac{n_1 \cdots n_{k}}{m_1 \cdots m_{k}} \right)^{it_1}
\left( \frac{n_1' \cdots n_{k}' }{m_1' \cdots m_{k}'  } \right)^{it_2 } \ dt_1 \ dt_2\\
&=N^{-2k}\int_{-T}^{T} \int_{-T}^{T} \sum_{{\bf n},{\bf m},{\bf n'},{\bf m'} \in \mathcal{Q}_k } 
\left( \frac{n_1 \cdots n_{k}}{m_1 \cdots m_{k}} \right)^{it_1 }  \left( \frac{n_1' \cdots n_{k}' }{m_1' \cdots m_{k}'} \right)^{it_2 } \ dt_1 \ dt_2,
\end{align*}
where, in the first line, the superscript $\sharp$ indicates a restriction to $\alpha$-regular variables $n_j, m_j,n_j',m_j'  \in \Gamma_{\alpha}(N)$. The summation in the last line runs over the set $\mathcal{Q}_k=\mathcal{Q}_k(N)$ consisting of those quadruples 
$({\bf n},{\bf m},{\bf n'},{\bf m'}) \in [1,N]^{4k}$ which are made up of  $\alpha$-regular components and satisfy the identity
$$  n_1 \cdots n_{k}  \cdot n_1' \cdots n_{k}' = m_1 \cdots m_{k} \cdot m_1' \cdots m_{k}'. $$
It follows that 
\begin{equation}
\Var[M_{k,\alpha}]=N^{-2k} \int_{-T}^{T} \int_{-T}^{T} \sum_{{\bf n},{\bf m},{\bf n'},{\bf m'} \in \mathcal{S}_k } 
\left( \frac{n_1 \cdots n_{k}}{m_1 \cdots m_{k}} \right)^{it_1 }  \left( \frac{n_1' \cdots n_{k}' }{m_1' \cdots m_{k}'} \right)^{it_2 } \ dt_1 \ dt_2,
\end{equation}
where $\mathcal{S}_k \subset \mathcal{Q}_k$ is made up of quadruples $({\bf n},{\bf m},{\bf n'},{\bf m'})$ for which $n_1 \cdots n_{k} \neq m_1 \cdots m_{k}$ 
(and hence $n_1' \cdots n_{k}' \neq m_1' \cdots m_{k}'$).\\

Next we separate $ \mathcal{S}_k $ into two parts. Let $\mathcal{S}_k^{-}$ contain those quadruples satisfying 
$|\log  \frac{n_1 \cdots n_{k}}{m_1 \cdots m_{k}}| \leq T^{-1/2} $
and write $\mathcal{S}_k^{+}$ for the complement of $\mathcal{S}_k^{-}$ inside $\mathcal{S}_k$. Accordingly, we write $\Var[M_k]=V^{-}+V^{+}$ to denote the resulting double integrals.\\
To treat $V^{+}$, we integrate with respect to $t_1$ (and treat the integration over $t_2$ trivially) to get that
\begin{align}\label{vplus}
N^{2k} |V^{+}|&= \bigg| \int_{-T}^{T} \int_{-T}^{T} \sum_{{\bf n},{\bf m},{\bf n'},{\bf m'} \in \mathcal{S}_k^{+}} 
\left( \frac{n_1 \cdots n_{k}}{m_1 \cdots m_{k}} \right)^{it_1 }  \left( \frac{n_1' \cdots n_{k}' }{m_1' \cdots m_{k}'} \right)^{it_2 } \ dt_1 \ dt_2 \bigg| \notag \\
& \ll T^{3/2}  \sum_{{\bf n},{\bf m},{\bf n'},{\bf m'} \in \mathcal{S}_k^{+}} 1 \ll  T^{3/2} \sum_{B \leq N^{2k}} \tau_{2k}(B)^2 \ll T^{3/2} N^{2k} (4 k \log N)^{4k^2 -1}.
\end{align}

Moving on to the treatment of $V^{-}$, it will be enough to give an upper bound for the cardinality of $\mathcal{S}_k^{-}$. In order to count the number of quadruples $({\bf n},{\bf m},{\bf n'},{\bf m'}) \in \mathcal{S}_k^{-}$, we may assume without loss of generality that 
$$d_1:=n_1 \cdots n_{k}<m_1 \cdots m_{k}=:d_2.$$
Since $T \geq 4$, we have that $d_2/d_1 \leq 1+2T^{-1/2}$ whenever $\log(d_2/d_1)\leq T^{-1/2}$ and as a result, we gather that
\begin{equation}\label{skminus}
|\mathcal{S}_k^{-}|  \leq \sum_{\substack{ 1 \leq d_1<d_2\leq N^k \\ d_2/d_1 \leq 1+2T^{-1/2}   } }^{\star} \  \sum_{\substack{ B \leq N^{2k}\\d_1,d_2 | B}} \tau_{2k,N; \alpha}(B)^2 =: \mathcal{T},
\end{equation}
where the starred sum is restricted to pairs $d_1, d_2 \in  \Gamma_{\alpha}^k(N)$. To deal with the expression $\mathcal{T}$ given just above, let us first extract the largest common divisor of $d_1$ and $d_2$. We write
$$d_1d_2=s^2d_1'd_2', \qquad s=\gcd(d_1,d_2), \qquad \gcd(d_1',d_2')=1.$$
Next, observe that for any pair of naturals $d_1'<d_2'$ satisfying $d_2'/d_1' \leq 1 +2T^{-1/2} $, we necessarily have that $d_1'\geq T^{1/2}/2$. Reordering the inner-most sum in $\mathcal{T}$, we first see that
\begin{equation*}
\sum_{\substack{ B \leq N^{2k}\\d_1,d_2 | B}} \tau_{2k,N; \alpha}(B)^2 =  
\sum_{\substack{ B \leq N^{2k}\\  s d_1'  d_2' | B}} \tau_{2k,N; \alpha}(B)^2  
\end{equation*}
and hence
\begin{equation}  \label{Tlastline}
\mathcal{T} \leq \sum_{s \leq N^k }   \ \sum_{\substack{ T^{1/2}/2 \leq  d_1' <d_2' \leq N^k/s  \\d_2'/d_1' \leq 1 +2T^{-1/2}  } }^{\star} \  
\sum_{ K \leq N^{2k}/s d_1'  d_2' } \tau_{2k,N; \alpha}(s d_1'  d_2'  K)^2,
\end{equation}
where we have once again restricted to variables $d_1', d_2' \in  \Gamma_{\alpha}^k(N)$.
To bound the triple sum in the last line, we first separate the variables $s,d_1',d_2'$ and $K$ by way of \eqref{basictaubounds} to find that 
$$\tau_{2k,N; \alpha}(s d_1'  d_2'  K)^2 \leq \tau_{4k^2}(s d_1'  d_2'  K) \leq \tau_{4k^2}(s) \tau_{4k^2}( d_1') \tau_{4k^2}( d_2') \tau_{4k^2}( K) $$

after which we estimate the sum over $K$ using \eqref{divisorpowersum}. This yields an `inner-most contribution'
$$\sum_{ K \leq N^{2k}/s d_1'  d_2' } \tau_{4k^2}( K)  \leq \frac{N^{2k} }{s d_1'  d_2' } (2 \log N^{2k})^{4k^2} \leq  \frac{N^{2k} }{s (d_1')^2  } (4 \log N^{k})^{4k^2} $$
to \eqref{Tlastline}. Next we observe that the variable $d_2'\leq N^k$ runs over integers for which $\Omega(d_2')\leq k (\log N)^{\alpha}\leq (\log N^k)^{\beta}$  (with some $\beta \in (0,1)$ depending on $\alpha$).  In other words, $d_2' \in \Gamma_{\beta}(N^k)$. Before applying Proposition \ref{mainuppershort} to the summation over $d_2'$, we recall that (see Remark \ref{sigmaremark}) the estimate \eqref{shortdivisor} is very poor when the summation interval is short and $k$ is large. For this reason it will be convenient to lengthen the range of $d_2'$ somewhat. Writing $\rho=(1000 \log_2 N)^{-1}$ and $\sigma=1-2\rho$, we may now combine Proposition \ref{mainuppershort} with a double application of \eqref{divisorpowersum} (over dyadic ranges) to find that
\begin{align*}
\mathcal{T}& \leq N^{2k} (4 \log N^k)^{4k^2} \sum_{s \leq N^k } \frac{\tau_{4k^2}(s) }{s}  \sum_{ T^{1/2} \leq d_1' \leq N^k } \frac{\tau_{4k^2}(d_1') }{(d_1')^2}\  
\sum_{\substack{ d_1'<d_2' \leq N^k \\ d_2' \in \Gamma_{\beta}(N^k) \\ d_2'/d_1'=1 +O(T^{-\rho } ) }} \tau_{4k^2}( d_2'  )\\
& \ll  N^{2k}(4  \log N^k)^{4k^2}\exp(20 (4k^2)^{1/\sigma^2} \log_2 N^k) T^{-\rho} \sum_{s \leq N^k }  \frac{\tau_{4k^2}(s) }{s}  \sum_{ d_1' \leq N^k } \frac{\tau_{4k^2}(d_1') }{d_1'} \\
&\ll N^{2k}(4 \log N^k)^{12k^2}\exp(81k^2 \log_2 N^k)  T^{-\rho}. 
\end{align*}

In summary, we have found that
$$|\mathcal{S}_k^{-}| \ll N^{2k}\exp(100k^2 \log_2 N^k) / T^{\rho}.$$
As an immediate consequence we see that 
\begin{align}\label{vminus}
N^{2k}|V^{-}|&= \bigg| \int_{-T}^{T} \int_{-T}^{T} \sum_{{\bf n},{\bf m},{\bf n'},{\bf m'} \in \mathcal{S}_k^{-}} 
\left( \frac{n_1 \cdots n_{k}}{m_1 \cdots m_{k}} \right)^{it_1 }  \left( \frac{n_1' \cdots n_{k}' }{m_1' \cdots m_{k}'} \right)^{it_2 } \ dt_1 \ dt_2 \bigg| \notag \\
& \ll T^{2-\rho}  N^{2k}  \exp(200k^2 \log_2 N) 
\end{align}
and hence, combining \eqref{vplus} and \eqref{vminus} we get that
\begin{align}
\Var[M_{k,\alpha}]&=V^{-}+V^{+} \ll T^{3/2} (4 k \log N)^{4k^2 -1}+  T^{2-\rho}   \exp(200k^2 \log_2 N)  \notag \\
& \ll T^{2-\rho}   \exp(200k^2 \log_2 N).
\end{align}

\end{proof}
\end{prop}

{\bf Concluding the proof of the lower bound in \eqref{mainestimate}. } Given $\gamma \in (0,1)$, we fix an $\alpha \in (\frac{1}{2}(1+ \gamma),1)$ and suppose that $C(x)$ satisfies \eqref{Cgrowth}.  We write 
$$ T=N^{C(N)}, \qquad \qquad \rho=\frac{1}{1000 \log_2 N}$$
and set $k= \lfloor (C(N) \log N)^{1/2}/(500 \log_2 N) \rfloor$. With this choice of parameters, the variance $\Var[M_{k, \alpha}]$ is well controlled since 
$$\Var[M_{k,\alpha}] \ll T^{2-\rho}    \exp(200k^2 \log_2 N) =o(T^2 ),$$
as $N \rightarrow \infty$, whereas \eqref{lowerdivisork} and \eqref{expectationmkalpha} give the lower bound
$$\E[M_{k,\alpha}] \gg T \exp \left( \frac{k^2}{200 \log_2 N } \right). $$
By Chebychev's inequality,
$$
\mathbb{P} \Big( |M_{k,\alpha}- \E[M_{k,\alpha}] | \,\geq\, \tfrac{1}{2} \E[M_k] \Big) \leq \frac{4 \Var[M_{k,\alpha}]}{\E[M_{k,\alpha}]^2} = o(1),
$$
from which it follows that $M_{k,\alpha} \geq\E[M_{k,\alpha}] /2$ with probability $1-o(1)$.  As a result we see that
$$
\sup_{|t|\leq T} |D_N^{\alpha}(t)|^{2k} \geq \frac{1}{2T} \int_{-T}^{T} |D_N^{\alpha}(t)|^{2k} \ dt =\frac{M_{k,\alpha}}{2T} \gg  \exp \left( \frac{k^2}{200 \log_2 N} \right),
$$
with probability $1-o(1)$ which, combined with \eqref{badsmall}, yields the desired lower bound \eqref{mainestimate}, with $B = 5 \cdot 10^{-6}$.

\subsection{Contrasting the independent variable case}\label{contrast}

In this final section we very briefly touch on the estimate \eqref{independentest} for independent random variables and explain why the sup-norm $\|D_N\|_{\infty}=\sup_{t \in \R} |D_N(t)|$
is larger than in the independent case. To address the second issue we must first recall the Bohr correspondence.\\ 
Letting $r=\pi(N)$, we define for each prime $p_j \leq N$ a complex variable $z_{p_j}$ on the unit circle $\T$ and write $\underline{z}=(z_{p_1},...,z_{p_r})$. We may then convert the Dirichlet polynomial $\mathcal{D}_N(s)=\sum_{n \leq N} X(n)n^{-s}$  into a trigonometric polynomial $Q(\underline{z})$ in $r$ variables as follows: replace each monomial $(p_{i_1}\cdots p_{i_j})^{-s}$ appearing in $\mathcal{D}_N(s)$ with the corresponding monomial $z_{p_{i_1}}\cdots z_{p_{i_j} }$. Under this identification one has Bohr's identity \cite[Eq (4.4.2.)]{QQ}
\begin{equation}\label{Bohrsup}
\|D_N\|_{\infty}=\sup_{\underline{z} \in \T^r} |Q(\underline{z})|.
\end{equation}     
Since the sequence $X(n)$ is completely multiplicative and takes values on the circle $\T$, the supremum on the RHS of \eqref{Bohrsup} occurs when $z_p=\overline{X(p)}$ at each prime $p$. Indeed,  in this case $\|D_N\|_{\infty} \geq |Q(\underline{z})|=  N$, which obviously matches the trivial upper bound $\|D_N\|_{\infty} \leq N$.

Moving on to the sup-norm estimate \eqref{independentest}, let $(r_n)_{n\in \N}$ denote a sequence of Steinhaus  i.i.d. variables. We observe that the moments of
\begin{equation*}
R_N(t)=\sum_{n \leq N} r_n n^{it}
\end{equation*}
are computed in much the same way as \eqref{expectationmk} and one obtains the straightforward evaluation 
\begin{align}\label{rnmoment}
\E \Big[ \int_{-T}^{T}| R_N(t)|^{2k} \ dt \Big] &=\int_{-T}^{T} \E \bigg[ \sum_{{\bf n},{\bf m} \in [1,N]^k } r_{n_1} \cdots r_{n_k} \overline{r_{m_1} \cdots r_{m_k} } \bigg]
\left( \frac{n_1 \cdots n_{k}}{m_1 \cdots m_{k}} \right)^{it } \ dt  \notag \\
& \sim2 k! T N^k
\end{align}
since only tuples $\left\{n_1,...,n_k \right\}$ which match up pairwise with the $\left\{m_1,...,m_k \right\}$ make a non-zero contribution in the above calculation. Combining the moment evaluation \eqref{rnmoment} with \eqref{sup}, we let $T=N^C$ and gather that
$$\E \Big[ \sup_{|t| \leq T} |R_N(t)| \Big]  \ll   \left(  k!    T N^{k+1} \log N   \right)^{1/2k}.$$
When $C$ is bounded away from zero, say $C \geq 1$, we choose an exponent $k \asymp C \log N$ and apply Chebychev's inequality to find that $\mathbb{P}(\sup_{|t| \leq N^C} |R_N(t)| \geq \lambda \sqrt{C N \log N})=O(1/\lambda)$, as claimed in 
\eqref{independentest}.

\bibliographystyle{siam}
\bibliography{RMFDirichlet} 
\end{document}